\documentclass [12pt] {report}
\usepackage{amssymb}
\newcommand {\D}[2] {\displaystyle\frac{\partial{#1}}{\partial{#2}}}

\newcommand {\Dd}[3] {\displaystyle\frac{\partial^2{#1}}{\partial{#2}\partial{#3}}}
\newcommand {\al} {\alpha}

\newcommand {\ga} {\gamma}
\newcommand {\la} {\lambda}

\newcommand {\de} {\delta}

\newcommand {\fr} {\displaystyle\frac}
\newcommand {\wt} {\widetilde}
\newcommand {\be} {\begin{equation}}
\newcommand {\ee} {\end{equation}}
\newcommand {\ba} {\begin{array}}
\newcommand {\ea} {\end{array}}
\newcommand {\bp} {\begin{picture}}
\newcommand {\ep} {\end{picture}}
\newcommand {\bc} {\begin{center}}
\newcommand {\ec} {\end{center}}
\newcommand {\bt} {\begin{tabular}}
\newcommand {\et} {\end{tabular}}
\newcommand {\lf} {\left}
\newcommand {\rg} {\right}

\newcommand {\cF} {{\cal F}}

\newcommand {\ses} {\medskip}

\newcommand {\bibit} {\bibitem}
\newcommand {\nin} {\noindent}

\newcommand {\De} {\Delta}

\usepackage{epsfig, latexsym,graphicx}
 \usepackage{amsmath}

\def\2#1#2#3{{#1}_{#2}\hspace{0pt}^{#3}}
\def\3#1#2#3#4{{#1}_{#2}\hspace{0pt}^{#3}\hspace{0pt}_{#4}}
\newcounter{sctn}
\def\sec#1.#2\par{\setcounter{sctn}{#1}\setcounter{equation}{0}
                  \noindent{\bf\boldmath#1.#2}\bigskip\par}

\begin {document}

\begin {titlepage}

 \ses\ses\ses\ses

\vspace{0.1in}

\begin{center}

{\large \bf FINSLEROID--FINSLER  PARALLELISM
}

\end{center}

\vspace{0.3in}

\begin{center}

\vspace{.15in} {\large G.S. ASANOV\\} \vspace{.25in}
{\it Division of Theoretical Physics, Moscow State University\\
119992 Moscow, Russia\\
{\rm (}e-mail: asanov@newmail.ru{\rm )}} \vspace{.05in}

\end{center}

\begin{abstract}

The Finsleroid--induced   scalar product, and hence the angle,  proves to remain unchanged under
the Finsleroid--type  parallel transportation
of involved vectors in the Landsberg case.
The two--vector extension of the Finsleroid metric tensor is proposed.

\ses

\nin {\it  Key words}: Finsler geometry, metric spaces, angle,
scalar product,  parallelism, relativity.

\end{abstract}

\end{titlepage}

\vskip 1cm

\vskip 1cm

\setcounter{sctn}{1} \setcounter{equation}{0}

{\large\bf 1. Introduction and synopsis}

\bigskip

The principal position of the Riemannian geometry is the phenomenon that  the angle between vectors does not change
under the parallel transportation of the vectors.
The theory of connection in the Riemannian geometry is developed, and taught to students,  subject to this observation.
 Can the phenomenon be transgressed from the Riemannian geometry to the Finsler geometry?
No transparent and constructive  answer is suggested
by the content of current literature
devoted to  Finsler spaces (see the books [1--3]).
This  notwithstanding, quite certain  positive answer proves to be a truth in the domain of the Finsleroid--Finsler
geometry  outlined in [4--8]. The answer is gained in the following succession  of steps. Firstly, we  use
the  Finsleroid--produced  scalar product ${\langle y_1,y_2\rangle_x}$ obtained on attentive studying
 the equations of geodesics in tangent spaces.
 Secondly,  we define the parallel
 displacement (1.1)
 of  vector $y^i$ with the respective spray--induced coefficients $\bar G^i{}_m$.
  Thirdly, we apply the Landsberg--case spray coefficients.
Lastly, we verify by straightforward
calculations (which are short and easy) that such a procedure does not change the value of
 ${\langle y_1,y_2\rangle_x}$  and, hence, the  Finsleroid angle.

Therefore, fixing the Landsberg case, we are entitled to conclude that
the Finsleroid approach proves to overcome the vague opinion that in
the Finsler geometry scientists may be ``in principle  equipped with only a family of
Minkowski norms'', so that ``yardsticks are assigned, but protractors are not''.
They can be equipped also with a convenient family of the two--vector products  ${\langle y_1,y_2\rangle_x}$,
 thereby with ``protractors''!

Sometimes the lack of two--vector angle is even  lifted ``to a high level of the principle of vintage''. The author of the
present paper (and not he alone) has
heard and read many times of ``the  specific grounds  that conclusively erase the concept of two--vector angle from the
 Finsler geometry, at least in the dimensions $N>2$'', and even of  the deep--wisdom advises
``better to forget of  two--vector angle when  opening the door to enter the
Finsler geometry!''. Secs. 1.6 and 1.7 of H.Rund's book [1] are tortured, much
and much, with ephemeral definitions of trigonometric functions and angles...

{

Let us try  only to be self--consistent!
 In the dimension $N=2$,
the  two--vector angle {\it does} enter the Finsler geometry, namely being  the
Landsberg angle. The obvious definition is the advantage of this
angle, with quite a similar significance as in the
two--dimensional Riemannian geometry, however in general
  there exists no  possibility to represent the angle
in the form of  an explicit algebraic function of two vectors. So, in the dimension $N=2$, the problem
with angle is of analytical, not conceptual,  nature, and we are to conclude that
 the two--dimensional Finsler geometry is a geometry!

Let us move in the dimensions $N > 2$.         The Finslerian indicatrix,
--- the extension of the Euclidean unit sphere, ---
 is  at our disposal.
 Nobody prevents
us from measuring angle between any two common--origin vectors by means of the
length
 of  the respective  arcs cut by vectors (or their continuations) from   unit circles
 located on the indicatrix,
  --- just in compliance with the known Euclidean school  methods.
 Again, the problem of getting the result may be only of  analytical nature.

\ses

 This circumstance thrusts  forth new
 questions fundamental
 to the very Realm of the nowaday Finsler geometry: should we
 consider that geometry ``old--fashioned'' from the new advantageous standpoint that is proposed by
the Finsleroid--induced  geometry? The vantage--ground answer is ``No'' in many principle aspects, particularly
the concepts of the Finslerian metric function and metric tensor, the Cartan tensor,
the geodesics and spray coefficients, the significance and geometry of indicatrix, the nonlinear covariant
derivative,
the connection and curvature on the tangent bundle, the flag curvature, etc.,  are keeping fine.
Simultaneously, the answer is decisively ``Yes'' in numerous new respects, including the occurrence
of the scalar product $\langle y_1,y_2\rangle_x$ between two vectors
 from which  many new  categories of the  cardinal geometrical
nature proper are stemming up.

{

Among such categories, ``the transportation preserving angle
between two vectors'' is notable, --- and can be consistently and  explicitly
tractable.
Indeed,
we may define  the {\it Finsleroid--Finsler  covariant differential}  $\de y$
of a vector $y$ along (horizontal) $dx$ in the natural way
\be
\de y^i~:=dy^i+\bar G^i{}_k(x,y)dx^k,
\ee
where
$
\bar G^i{}_k=\frac 12 \partial  G^i/\partial y^k
$
and
$
G^i~:=\ga^i{}_{mn}y^my^n$
are the respective spray coefficients,
with
 $\ga^i{}_{mn}$
standing for  the Finslerian Christoffel symbols constructed from the Finsleroid--Finsler  metric function $K$.
The vector $y$ is said to undergo the {\it parallel transportation along  $dx$} if
\be
\de y=0.
\ee

Then
analytically the condition for the scalar product
 $\langle y_1,y_2\rangle_x$
 to be unchanged under such a transportation of vectors $y_1,y_2$
reads
\be
\D{{\langle y_1,y_2\rangle}_x}{x^k}
-\bar G^n{}_k(x,y_1)\D{{\langle y_1,y_2\rangle}_x}{y_1^n}
-\bar G^n{}_k(x,y_2)\D{{\langle y_1,y_2\rangle}_x}{y_2^n}=0.
\ee

We claim

\ses\ses

{\large Parallel Transportation  Theorem.} In the Landsberg case of the Finsleroid--Finsler space,
the condition (1.3) holds fine.

\ses\ses

The proof is arrived at after direct  calculations which are not lengthy, as will be demonstrated in
Appendix A.
Thus,  both the scalar product (angle) of pair of vectors as well the parallel transportation
retaining the product (angle) can nicely be transgressed
from the Riemannian geometry to the Finsleroid--Finsler
geometry in a simple analytical way.
{\it This parallelism in the Finsleroid domain comes to play  replacing the Levi-Civita parallelism functioned
conventionally  in the Riemannian geometry.}

In Finsler geometry,
we have two concepts of vector length.
Namely,
 we can use the Finslerian metric function $F = \sqrt{g_{ij}(x, y)y^iy^j}$ to assign the {\it absolute length}
\be
||y||_x= \sqrt{g_{ij}(x, y)y^iy^j}
\ee
to the vector $y\in T_xM$.
Simultaneously, taken another vector $\wt y$ in the same tangent space, such that $\wt y,y\in T_xM$,
the Finsler geometry theory [1] provides us with  the {\it relative length}
\be
||y||_{\wt y}=\sqrt{g_{ij}(x,\wt y)y^iy^j}
\ee
which measures the vector $y$ relative to a supporting vector $\wt y$.

{

It is natural to wonder whether  the length definitions (1.4) and (1.5) can be extended to give us
 respective  notions of angles.
The second case is extending in quite an obvious and traditional way as follows:
at any fixed point $x$, we have  the {\it relative  scalar product}
\be
\langle y_1,y_2\rangle_y=g_{ij}(x,y)y_1^iy_2^j  \qquad  \text{in  any Finsler  space}.
\ee
No possibility to extend properly  the absolute case (1.4)  is proposed in the books [1--3]
(and in the current literature).
However, the Finsleroid--Finsler geometry is wonderful in that it provides us with the following
{\it absolute  scalar product}:
\be
\langle y_1,y_2\rangle_x=G_{ij}(x,y_1,y_2)y_1^iy_2^j \qquad  \text{in  the Finsleroid--Finsler  space}
\ee
(see (2.6)).

The occurrence of the scalar product $\langle y_1,y_2\rangle_x$
suggests naturally proposing  two--vector extensions $Y_{1i}(x,y_1,y_2),\,Y_{2i}(x,y_1,y_2),\,G_{ij}(x,y_1,y_2)$,
$A_{a;ijk}(x,y_1,y_2),a=1,2,$ of the ordinary Finslerian definitions of the covariant vector
 $y_i=K\partial K/\partial {y^i}$,
 the  metric tensor $g_{ij}=g_{ij}(x,y)$, and
 the  Cartan tensor $A_{ijk}(x,y)$.
The explicit components of these extensions are found in Section 2. By their use,  direct calculations
(shown in Appendix A) reveal the validity of the following theorem.

\ses\ses

{\large Match  Theorem.} In the Finsleroid--Finsler space under study, the  following limits are fulfilled:
\be
\lim_{y_2\to y_1=y}{Y_{1i}}=
\lim_{y_2\to y_1=y}{Y_{2i}}= y_i
\ee
and
\be
\lim_{y_2\to y_1=y}{G_{ij}}=g_{ij},
\ee
together with
\be
\lim_{y_2\to y_1=y}{A_{a;ijk}}=A_{ijk}.
\ee

\ses\ses


Therefore, we are to expect that the Finsleroid geometry theory is not of
primarily complete nature and should be regarded as a limiting $(y_1\,=\,y_2)$--case  of the respective
two--vector extended theory to be developed in future.
The important  nature of the tensor $G_{ij}$ can be seen in the equality (1.7) which expresses the scalar product
$\langle y_1,y_2\rangle_x$ by means of the tensor.

\ses

Whether the lengths and scalar products (1.4)--(1.7)
 remain unchanged under the parallel transportation of the involved vectors?
 Due answers
 will be formulated in Section 3, yielding
 {\bf the Total Category of Parallelism}, in which all the distant--parallelism concepts
are meaningful as well as  representable in an explicit and simple way, applying the Landsberg case.

 {

Throughout the paper, the notation is the same as in the previous work [4--8],
in which
we have introduced the Finsleroid--Finsler space $\cF\cF^{PD}_{g}$
under the  condition
that  the norm $||b||$ of  the Finsleroid--axis 1-form
\be
b=b_iy^i
\ee
 is equal to 1:
 \be
 a_{ij}(x)b^i(x)b^j(x)=1,
 \ee
 where $a_{ij}$ stands for the metric tensor of the associated Riemannian space metricized by
 the function $S(x,y)=\sqrt{ a_{ij}(x)y^iy^j}$ of point $x$ and tangent vector $y$.
We shall normalize the fundamental Finsleroid--Finsler metric function $K$ to  fulfill the  condition
\be
g_{ij}\bigl(x,b(x)\bigr)=a_{ij}(x)
\ee
which in turn entails
\be
K\bigl(x,b(x)\bigr)=1.
\ee
Thus, the Finsleroid--geometry properties come to play when the tangent vector $y^i$
begins deviating from the vector $b^i$.
The conditions (1.13) and (1.14)  assign actually the {\it  correspondence principle} to make comparison between
the Finsleroid--Finsler space and  the associated Riemannian space.
We have also
\be
g_{ij}(x,y)\bigl|_{g=0}\bigr.= a_{ij}(x), \qquad
K(x,y)\bigl|_{g=0}\bigr.= S(x,y).
\ee
The equalities (1.12)--(1.15) are essential in developing our subject.

\ses

When proceeding in this direction, should the recent Finslerian
theory of
 connection and curvature   be recapitulated anew to comply strictly with the  Levi--Civita idea?
May the nonlinear methods of construction of covariant
derivatives  start  coming to play significantly? All these questions are important
and open  to make deep inquiry.

The limitation of  our  parallel transportation theorem
is that we use the  spray coefficients $\bar G^i$ of the Landsberg
case,  --- and the present author has  not succeeded as yet in
answering the troublesome question whether the conclusion can be
extended to more general cases; some auxiliary calculations are
presented in  Appendix B.

The present paper deals everywhere with the positive--definite case. However, all the conclusions made
can directly be re--formulated to apply to the relativistic pseudo--Finsleroid--Finsler space.

{
\ses\ses

\setcounter{sctn}{2}
\setcounter{equation}{0}

{\large\bf 2. Scalar product and two--vector  tensors}

\vspace{0.7cm}

Below, we use a pair $y_1,y_2\in T_xM$ of tangent vectors  supported by a fixed point $x\in M$ of the  background
$N$--dimensional manifold $M$.

If a Finsler space involves a
 scalar product $\langle y_1,y_2\rangle_x$
which  possesses the homogeneity
\be
\langle ky_1,y_2\rangle_x=k\langle y_1,y_2\rangle_x,
\quad
\langle y_1,ky_2\rangle_x=k\langle y_1,y_2\rangle_x,
\qquad k>0,\, \forall y_1,y_2,
\ee
then it is attractive to explicate
the {\it two--vector covariant vectors}
\be
Y_{1i}(x,y_1,y_2)~:=\D{\langle y_1,y_2\rangle_x}{y_1^i},
\qquad
Y_{2j}(x,y_1,y_2)~:=\D{\langle y_1,y_2\rangle_x}{y_2^j}
\ee
and
\ses
the {\it two--vector metric tensor}
\be
G_{ij}(x,y_1,y_2)~:=\Dd{\langle y_1,y_2\rangle_x}{y_1^i}{y_2^j},
\ee
together with
 the following
{\it  two--vector  extension of the Cartan tensor}:
\be
A_{a;ijk}(x,y_1,y_2)=\langle y_1,y_2\rangle_xC_{a;ijk}(x,y_1,y_2), \qquad a=1,2,
\ee
with
\be
C_{1;kij}(x,y_1,y_2)~:=\D{G_{ij}(x,y_1,y_2)}{y_1^k},
\qquad
C_{2;ijk}(x,y_1,y_2)~:=\D{G_{ij}(x,y_1,y_2)}{y_2^k}.
\ee

The homogeneity (2.1)
entails obviously the  identities
\be
G_{ij}(x,y_1,y_2)y^i_1y^j_2=
y_1^iY_{1i}(x,y_1,y_2)=Y_{2i}(x,y_1,y_2)y_2^i=\langle y_1,y_2\rangle_x
\ee
and
\be
y_1^iG_{ij}(x,y_1,y_2)=Y_{2j}(x,y_1,y_2), \qquad G_{ij}(x,y_1,y_2)y_2^j=Y_{1i}(x,y_1,y_2).
\ee
The generalized symmetry
\be
G_{ij}(x,y_1,y_2)=G_{ji}(x,y_2,y_1)
\ee
is valid.

Also,
\be
y_1^kC_{1;kij}(x,y_1,y_2)=
C_{2;ijk}(x,y_1,y_2)y_2^k=0.
\ee

{

We shall mark  quantities by the subscript `1' if they are taken at the  value $y=y_1$,
{\it resp.}  by the subscript `2'  at the  value $y=y_2$, as exemplified by
\be
A_1\!=\!A(x,y_1), A_2\!=\!A(x,y_2), \,B_1\!=\!B(x,y_1), B_2\!=\!B(x,y_2),\, K_1\!=\!K(x,y_1), K_2\!=\!K(x,y_2).
\ee
The Finsleroid--Finsler scalar product was presented explicitly by the formulas (2.33) and (2.34) in [7],
such that
\be
\langle y_1,y_2\rangle_x=K_1K_2\cos\al_x
\ee
 and
\be
\al_x=\fr1h\arccos\la,
\ee
where
\be
\la= \fr{ A_1A_2+h^2r_{ij}y_1^iy_2^j} {\sqrt{B_1}\,\sqrt{B_2}}.
\ee
Introducing  the notation
\be
\ga=\fr{\sin\al}{\sqrt{1-{\la}^2}}\equiv \fr{\sin\al}{\sin(h\al)}
\ee
(and avoiding indication of the subscript $x$ for $\al$),
we deduce the explicit components
\be
Y_{1i}=
\D{K_1}{y_1^i}K_2\cos\al
+\fr{\ga}hK_1K_2
\D{\la}{y_1^i},
\qquad
Y_{2i}=
K_1\D{K_2}{y_2^i}\cos\al
+\fr{\ga}hK_1K_2\D{\la}{y_2^i},
\ee
and
\be
G_{ij}=    \fr1{K_1}      \D{K_1}{y_1^i}Y_{2j}
+\fr{\ga}hK_1   \D{K_2}{y_2^i}
\D{\la}{y_1^i}
+
\fr{\ga}hK_1K_2
\Dd{\la}{y_1^i}{y_2^j}
+\fr1hK_1K_2\D{\ga}{\la}
\D{\la}{y_1^i}\D{\la}{y_2^j},
\ee
or
\be
G_{ij}\!=\!    \fr1{K_1}      \D{K_1}{y_1^i}Y_{2j}    +     \fr1{K_2}      \D{K_2}{y_2^i}Y_{1j}
-
\D{K_1}{y_1^i}\D{K_2}{y_2^j}\cos\al
+
\fr{\ga}hK_1K_2
\Dd{\la}{y_1^i}{y_2^j}
+
\fr{1}hK_1K_2\D{\ga}{\la}
\D{\la}{y_1^i}\D{\la}{y_2^j}.
\ee

{

\ses\ses

\setcounter{sctn}{3}
\setcounter{equation}{0}

{\large\bf 3. Parallel transportation}

\vspace{0.7cm}

\ses

Let us  introduce   the parallel transportation of a vector $X\in T_xM$ along an infinitesimal (horizontal)
displacement $dx$ by following the known method described in Section 6.4 of [1].
Below, all the components $g_{ij}$
and
$
\bar G^h,  \bar G^h{}_{i}, \bar G^h{}_{ij}, \bar G^h{}_{kij}$
are implied to depend on the argument $(x,X)$.
The notation $\de X$ features
the  covariant differential (1.1), so that
\be
\de X^h=dX^h+\bar G^h{}_kdx^k
\equiv dX^h+\bar G^h{}_{km}X^mdx^k.
\ee
For the  covariant vector
\be
X_h=g_{hk}X^k
\ee
we take
\be
\de X_h=dX_h-X_l\bar G^l{}_{kh}dx^k,
\ee
such that
\be
\de (X_hX^h)=d(X_hX^h).
\ee

Proceeding in this way,
we introduce the {\it covariant differential of the Finslerian metric tensor}
\be
\de g_{ij}=\D{ g_{ij}}{x^k}dx^k
+\D{ g_{ij}}{X^h}dX^h
-g_{ih}\bar G^h{}_{jk}dx^k
-g_{jh}\bar G^h{}_{ik}dx^k,
\ee
which can also be written as
\be
\de g_{ij}=\fr{\de g_{ij}}{\de x^k}dx^k
+\D{ g_{ij}}{X^h}\de X^h
\ee
with
\be
\fr{\de g_{ij}}{\de x^k}=\D{g_{ij}}{x^k}
-\D{ g_{ij}}{X^h}\bar G^h{}_k
-g_{ih}\bar G^h{}_{jk}
-g_{jh}\bar G^h{}_{ik}.
\ee
Here the right--hand part can be transformed to yield
\be
\fr{\de g_{ij}}{\de x^k}=
X_h\bar G^h{}_{kij},
\ee
so that
\be
\de g_{ij}=
\D{ g_{ij}}{ X^h}\de X^h
+X_h\bar G^h{}_{kij}dx^k.
\ee
It is well--known that
\be
 \bar G^h{}_{kij}X^i=0 \quad \text{in any Finsler space}
 \ee
and
\be
 X_h\bar G^h{}_{kij}=0 \quad \text{in the Landsberg case of  Finsler space}.
 \ee
{



Accordingly, we introduce

\ses

DEFINITION. A vector $X$ is said to be {\it parallel} under the  displacement,
if $\de X=0$.
Also, the metric tensor
$g_{ij}=0$
{\it behaves parallel}, if $\de g_{ij}=0$ when
$\de X=0$.

\ses

\ses

NOTE. The equalities (3.7) and (3.8) are well--known from the book [1],
in which they were discussed as ``Berwald covariant derivative
of the Finslerian metric tensor", with the coefficients $\bar G^h{}_{ik}$ being treated as
 the ``Berwald connection coefficients" (see (3.10) of Section 3.3 in [1];
 the coefficients were denoted in [1] to read simply  $ G^h{}_{ik}$).
In our case,
the above formulas (3.9) and (3.5) are tantamount to, respectively,     Eqs. (4.17) and (4.16) of Section 6.4 of [1];
that Section was devoted to the nonlinear connection, so that we may qualify (3.1) by the status of
 the {\it nonlinear covariant differential of vector}.

\ses

By comparing (3.9) and (3.10) with (1.4) and (1.6), we just conclude that
\be
\de||y||_x=0 \quad \text{ under} \quad \de y= 0,
 \quad \text{in any Finsler space}
\ee
and
\be
\langle y_1,y\rangle_y=0 \quad \text{ under} \quad \de y= \de y_1=0,
 \quad \text{in any Finsler space},
\ee
where
\be
\langle y_1,y\rangle_y=g_{ij}(x,y)y_1^iy^j
\equiv y_1^iy_i.
\ee

However, the assertions of the type (3.12) and (3.13) are not applicable to the full scalar products.
The reason is that  the products involve  the Finslerian metric tensor
$g_{ij}(x,y)$
which, in contrast to
 the Riemannian metric tensor proper, depends on the transported vector $y$.
The parallelism property may be
a truth in the particular case when
the tensor $g_{ij}(x,y)$ itself is unchanged under the parallel transportation of the argument vector $y$.
Let  a set of vectors $\wt y, y,y_1,y_2$  be supported by same point $x$.
In view of the nullification (3.11),  the property said
 occurs as follows:
\be
\de g_{ij}=0 \quad \text{ under} \quad \de y=0 ,\quad \text{in the Landsberg case of Finsler space}.
\ee

{

This directly entails the assertion
\be
||y||_{\wt y}=0 ~~ \text{if} ~
~ \de \wt y= \de y=0, \quad
 \text{in the Landsberg case of Finsler  space}.
\ee
This chain is continuing as follows:
\be
\de\langle y_1,y_2\rangle_y=0 \quad \text{ under} \quad \de y= \de y_1= \de y_2=0,
 \quad \text{in the Landsberg case of Finsler space},
\ee
and
\be
\de\langle y_1,y_2\rangle_x=0 ~~ \text{if} ~
~ \de y_1= \de y_2=0, ~~
 \text{in the Landsberg case of the Finsleroid--Finsler  space}.
\ee
If we consider  the {\it relative angle}
\be
\al_y(y_1,y_2)=
\arccos
\fr{\langle y_1,y_2\rangle_y}{||y_1||_y||y_2||_y}  \qquad  \text{in  any Finsler  space}
\ee
and
the {\it absolute  angle}
\be
\al_x(y_1,y_2)=
\arccos
\fr{\langle y_1,y_2\rangle_x}{||y_1||_x||y_2||_x} \qquad  \text{in  the Finsleroid--Finsler  space},
\ee
from the above we are entitled to conclude that
\be
\de\al_y(y_1,y_2)=0 \quad \text{ under} \quad \de y= \de y_1= \de y_2=0,
 \quad \text{in the Landsberg case of Finsler space},
\ee
and
\be
\de\al_x(y_1,y_2)=0 ~~ \text{if} ~
~ \de y_1= \de y_2=0, ~~
 \text{in the Landsberg case of the Finsleroid--Finsler  space}.
\ee
Assuming the Landsberg case is essential.

Under  conditions of the previous assertion (3.22), the two--vector objects (2.2) and (2.3) are also parallel:
\be
\de G_{ij}=0
\ee
and
\be
\de Y_{1i}=\de Y_{2i}=0.
\ee

{

\ses\ses

\setcounter{equation}{0}

\bc
{\large \bf Appendix A. ~ Two--vector limits and parallelism condition}
\ec

\ses\ses

Let us  apply   the formulas
(2.11)--(2.13) to verify the Match Theorem of Section 1.
Using the derivatives
\be
\D b{b_i}=y^i, \qquad \D q{b_i}=-\fr bq y^i, \qquad \D A{b_i}=\fr{q-\fr12gb}qy^i, \qquad \D B{b_i}=g\fr{q^2-b^2}qy^i,
\ee
we obtain
\be
\D{\la}{b_i}=c_1y_1^i+c_2y_2^i
\ee
with
the coefficient
\be
c_1=
\fr1{q_1}
\Biggl[
\fr{A_2\lf(q_1-\fr12gb_1\rg)-h^2q_1b_2}{\sqrt{B_1}\,\sqrt{B_2}}
-\fr{g\la}{2B_1}(q_1^2-b_1^2)
\Biggr]
\ee
which can be simplified to read
\be
c_1=
\fr g{2q_1}
\Biggl[
\fr{q_1q_2-b_1b_2+\fr12g(q_1b_2-q_2b_1)}{\sqrt{B_1}\,\sqrt{B_2}}
-\fr{\la}{B_1}(q_1^2-b_1^2)
\Biggr].
\ee
The quantity $c_2$ is obtainable from $c_1$ by performing the subscript interchange
$1\leftrightarrow 2$.

{

\nin
With the help of the   notation
\be
t_{Ai}=a_{in}y_A^n+\fr12g(q_Ab_i+\frac{b_A}{q_A}v_{Ai}),\quad
 b_{Ai}=b_i+\fr 12g\fr{v_{Ai}}{q_A},
 \qquad A=1,2,
\ee
we find that
\be
\D{\la}{y_1^i}=
\fr{\lf(b_i+\fr12g\fr{v_{1i}}{q_1}\rg)A_2+h^2v_{2i}} {\sqrt{B_1}\,\sqrt{B_2}}
-\fr{\la t_{1i}}{B_1}, \quad
\D{\la}{y_2^i}=
\fr{\lf(b_i+\fr12g\fr{v_{2i}}{q_2}\rg)A_1+h^2v_{1i}} {\sqrt{B_1}\,\sqrt{B_2}}
-\fr{\la t_{2i}}{B_2}.
\ee
Contractions show that
\be
\D{\la}{y_1^i}y_1^i=
\D{\la}{y_2^i}y_2^i=0
\ee
\ses\ses
and
\be
\D{\la}{y_1^i}b^i=
\fr{A_2}{\sqrt{B_1}\,\sqrt{B_2}}
-\fr{\la A_1}{B_1},
\qquad
\D{\la}{y_2^i}b^i=
\fr{A_1}{\sqrt{B_1}\,\sqrt{B_2}}
-\fr{\la A_2}{B_2}.
\ee
Appropriate differentiation  yields
\be
\Dd{\la}{y_1^i}{y_2^j}=
\fr{b_{1i}b_{2j}
+h^2r_{ij}} {\sqrt{B_1}\,\sqrt{B_2}}
-\fr{b_{1i}A_2+h^2v_{2i}} {B_2\sqrt{B_1}\,\sqrt{B_2}}
t_{2j}
-\fr{1}{B_1}
t_{1i}
\Biggl(
\fr{b_{2j}A_1+h^2v_{1j}} {\sqrt{B_1}\,\sqrt{B_2}}
-\fr{\la}{B_2}
t_{2j}
\Biggr).
\ee
We may observe the properties
\be
\D{\la}{y_1^i}{\Bigl|_{y_2=y_1}\Bigr.}=\D{\la}{y_2^i}{\Bigl|_{y_2=y_1}\Bigr.}=0
\ee
and
\be
\Dd{\la}{y_1^i}{y_2^j}{\Bigl|\Bigr._{y_2=y_1=y}}=
\fr{\lf(b_i+\fr12g\fr{v_{i}}{q}\rg)
\lf(b_j+\fr12g\fr{v_{j}}{q}\rg)
+h^2r_{ij}} {B}
-\fr{t_it_j}{B^2},
\ee
\ses\\
where
$t_i=a_{in}y^n+\frac12g(qb_i+\frac{b}{q}v_i).$
If we compare the right--hand part of   (A.11) with the structure of the Finsleroid angular metric tensor $h_{ij}$
(see [5,7]), we obtain the simple equality
\be
\Dd{\la}{y_1^i}{y_2^j}{\Bigl|\Bigr._{y_2=y_1=y}}=\fr{h^2}{K^2}h_{ij}.
\ee
Taking into account (A.11) and the nullifications (A.7), together with the limits
\be
\lim_{\la\to 1} \fr{\sin\al}{\sqrt{1-{\la}^2}}=\fr1h,
\qquad
\lim_{\la\to 1}\D{\fr{\sin\al}{\sqrt{1-{\la}^2}}}{\la}=\fr{1-h^2}{h^3},
\ee
we are entitled to conclude from  the  formulas (2.15)--(2.17) that the claimed limits (1.8)--(1.10) of the theorem
are valid.

{

Now we turn to the Parallel Transportation  Theorem of Section 1.
Let us take two vectors $y_1,y_2\in T_xM$.
To establish  the vanishing (1.3),   we must apply accurate calculations to verify that
\be
\D{{\la}}{x^k}
-\fr12 G^n{}_k(x,y_1)\D{{\la}}{y_1^n}
-\fr12 G^n{}_k(x,y_2)\D{{\la}}{y_2^n}=0
\ee
with the function $\la$ given by (2.13), and with the Landsberg--case spray--induced coefficients
\be
 G^i{}_k=
\fr {gk}{q}
\biggl[
(u_k-bb_k)v^i+q^2(\de_k{}^i-b_kb^i)
\biggr]+2a^i{}_{km}y^m
\ee
(these coefficients can be found  in [5--8]).
Denoting
$$
G_1{}^i{}_k=G^i{}_k(x,y_1), \qquad  G_2{}^i{}_k=G^i{}_k(x,y_2),
$$
we obtain
\be
\D{\la}{y_1^i}G_1{}^i{}_k
=\fr {gk}{q_1}\D{\la}{y_1^i}\biggl[(q_1)^2\de_k{}^i
-b^im_{1k}
\biggr]+\De
=
gkq_1\D{\la}{y_1^k}
-\fr {gk}{q_1}
m_{1k}
\D{\la}{y_1^i}b^i+\De,
\ee
where $m_{1k}=b_1v_{1k}+(q_1)^2b_k$
and
 $\De$ symbolizes the summary of the terms which involve partial derivatives
of the input Riemannian metric tensor $a_{ij}$ with respect to the coordinate variables $x^k$.
 On  simplifying and applying  (A.6),  the right--hand part in (A.16) becomes
$$
gk
\biggl[
q_1\D{\la}{y_1^k}
+\fr {1}{q_1}m_{1k}
\fr{\la A_1}{B_1}
-\fr {1}{q_1}m_{1k}
\fr{A_2}{\sqrt{B_1}\,\sqrt{B_2}}
\biggr]+\De
$$
$$
=
gk
\Biggl[
q_1
\fr{\lf(b_k+\fr12g\fr{v_{1k}}{q_1}\rg)A_2+h^2v_{2k}} {\sqrt{B_1}\,\sqrt{B_2}}
-\fr{q_1\la}{B_1}t_{1k}
+\fr {1}{q_1}m_{1k}
\fr{\la A_1}{B_1}
-\fr {1}{q_1}m_{1k}
\fr{A_2}{\sqrt{B_1}\,\sqrt{B_2}}
\Biggr]
+\De.
$$
Here, all the terms proportional to $b_k$ are cancelled, leaving us with
$$
gk
\Biggl[
q_1
\fr{\fr12g\fr{v_{1k}}{q_1}A_2+h^2v_{2k}} {\sqrt{B_1}\,\sqrt{B_2}}
-\fr{q_1\la}{B_1}
\lf(1+
\fr12g\fr{b_1}{q_1}\rg)v_{1k}
+\fr {1}{q_1}b_1v_{1k}
\fr{\la A_1}{B_1}
-\fr {1}{q_1}b_1v_{1k}
\fr{A_2}{\sqrt{B_1}\,\sqrt{B_2}}
\Biggr]
+\De.
$$
{
Eventually,
\be
\D{\la}{y_1^i}G_1{}^i{}_k=
\fr{gk}{q_1}
\Biggl[
\fr{\fr12gq_1A_2v_{1k}+h^2(q_1)^2v_{2k}} {\sqrt{B_1}\,\sqrt{B_2}}
-\fr{q_1\la(q_1+\fr12gb_1)v_{1k}}{B_1}
+
\fr{\la A_1b_1v_{1k}}{B_1}
-
\fr{A_2b_1v_{1k}}{\sqrt{B_1}\,\sqrt{B_2}}
\Biggr]+\De.
\ee
Interchanging here
 ${\scriptstyle 1\leftrightarrow2}$
 yields the quantity
$\frac{\partial \la}{\partial y_2^i}{\scriptstyle G_2{}^i{}_k}$.
Now, using the characteristic Landsberg condition
$
\nabla_ib_j=k(a_{ij}-b_ib_j),
$
we get
\be
\D{\la}{x^k}=k(c_1v_{1k}+c_2v_{2k})+\De,
\ee
where (A.4) should be used.
With the formulas (A.17) and (A.18), the validity of the vanishing (A.14) can readily be seen.

{

\ses\ses

\setcounter{equation}{0}

\bc
{\large \bf Appendix B. ~ Use of full spray coefficients}
\ec

\ses\ses

Suppressing  the Landsberg condition, the full spray coefficients are given by
the representation (A.48) of [7]  which yields
\ses\\
\be
G_1{}^i{}_k=
gP_{1k}v^i_1+gQ_1(\de^i{}_k-b^ib_{k})
-g\fr{v_{1k}}{q_1}f_1^i-gq_1f^i{}_k+2a^i{}_{km}y_1^m,
\ee
where
\be
P_{1k}=
-\fr1{q_1^3}v_{1k}  y_1^jy_1^h\nabla_jb_h
+\fr1{q_1} y_1^j(\nabla_jb_k+\nabla_kb_j)
+gb^j\nabla_jb_k
\ee
and
\be
Q_1=
\fr1{q_1} y_1^jy_1^h\nabla_jb_h
+gy_1^hb^j\nabla_jb_h.
\ee
We use the  notation
\be
f^i=f^i{}_ny^n,\qquad
f^i{}_n=a^{ik}f_{kn}, \qquad
f_{mn}=
\nabla_mb_n-\nabla_nb_m
\equiv \D{ b_n}{x^m}-\D {b_m}{x^n},
\ee
where
 the nabla  means the covariant derivative in terms of the associated Riemannian space.
We obtain
\be
\D{\la}{y_1^i}G_1{}^i{}_k
=gP_{1k}v^i_1\D{\la}{y_1^i}
+gQ_1(\de^i{}_k-b^ib_{k})\D{\la}{y_1^i}
-g\fr{v_{1k}}{q_1}f_1^i
\D{\la}{y_1^i}
-gq_1f^i{}_k
\D{\la}{y_1^i}
+\De,
\ee
or
\ses\\
\be
\D{\la}{y_1^i}G_1{}^i{}_k
=gQ_1\D{\la}{y_1^k}
+gP_{1k}\D{\la}{y_1^i}v^i_1
-gQ_1b_{1k}\D{\la}{y_1^i}b^i
-g\fr{v_{1k}}{q_1}
\D{\la}{y_1^i}f_1^i
-gq_1
\D{\la}{y_1^i}f^i{}_k
+\De.
\ee
The representation (B.6) extends the formula (A.16) of the preceding Appendix A.

{

Now we start calculating in the straightforward way:
\be
\D{\la}{y_1^i}G_1{}^i{}_k
=gQ_1\D{\la}{y_1^k}
+gP_{1k}\D{\la}{y_1^i}v^i_1
-gQ_1b_{k}\D{\la}{y_1^i}b^i
-g\fr{v_{1k}}{q_1}
\D{\la}{y_1^i}f_1^i
-gq_1
\D{\la}{y_1^i}f^i{}_k
+\De,
\ee
or after required insertions
\ses\\
$$
\D{\la}{y_1^i}G_1{}^i{}_k
=
gQ_1
\Biggl(
\fr{b_{1k}A_2+h^2v_{2k}} {\sqrt{B_1}\,\sqrt{B_2}}
-\fr{\la}{B_1}
t_{1k}
\Biggr)
+
gP_{1k}
\Biggl(
\fr{b_{1i}A_2+h^2v_{2i}} {\sqrt{B_1}\,\sqrt{B_2}}
-\fr{\la}{B_1}
t_{1i}
\Biggr)
v_1^i
$$
\ses
$$
-gQ_1b_{1k}
\Biggl(
\fr{b_{1i}A_2+h^2v_{2i}} {\sqrt{B_1}\,\sqrt{B_2}}
-\fr{\la}{B_1}
t_{1i}
\Biggr)
b^i
-g\fr{v_{1k}}{q_1}
\Biggl(
\fr{b_{1i}A_2+h^2v_{2i}} {\sqrt{B_1}\,\sqrt{B_2}}
-\fr{\la}{B_1}
t_{1i}
\Biggr)
f_1^i
$$
\ses
$$
-gq_1
\Biggl(
\fr{b_{1i}A_2+h^2v_{2i}} {\sqrt{B_1}\,\sqrt{B_2}}
-\fr{\la}{B_1}
t_{1i}
\Biggr)
f^i{}_k
+\De,
$$
where the notation (A.5) has been applied.
Simplifying yields
\ses\\
$$
\D{\la}{y_1^i}G_1{}^i{}_k
=
gQ_1
\Biggl\{
\fr{\fr12g\fr{v_{1k}}{q_1}A_2+h^2v_{2k}} {\sqrt{B_1}\,\sqrt{B_2}}
-\fr{\la}{B_1}
\lf(1+\fr12g\fr{b_1}{q_1}\rg)
v_{1k}
\Biggr\}
$$
\ses
$$
+
gP_{1k}
\Biggl\{
\fr{\fr12g\fr{v_{1i}}{q_1}A_2+h^2v_{2i}} {\sqrt{B_1}\,\sqrt{B_2}}
-\fr{\la}{B_1}
\biggl[
v_{1i}+
\fr12g\fr{b_1}{q_1}v_{1i}
\biggr]
\Biggr\}
v_1^i
$$
\ses
\ses
$$
-g\fr{v_{1k}}{q_1}
\Biggl\{
\fr{\lf(b_i-\fr12g\fr{b_1b_i}{q_1}\rg)A_2+h^2v_{2i}} {\sqrt{B_1}\,\sqrt{B_2}}
-\fr{\la}{B_1}
\fr12g\Bigl(q_1b_i-\fr{b_1}{q_1}b_1b_i\Bigr)
\biggr]
\Biggr\}
f_1^i
$$
\ses
$$
-gq_1
\Biggl\{
\fr{\lf(b_i+\fr12g\fr{v_{1i}}{q_1}\rg)A_2+h^2v_{2i}} {\sqrt{B_1}\,\sqrt{B_2}}
-\fr{\la}{B_1}
\biggl[
u_{1i}+
\fr12g\Bigl(q_1b_i+\fr{b_1}{q_1}v_{1i}\Bigr)
\biggr]
\Biggr\}f^i{}_k
+\De,
$$
{
\nin
or
\ses\\
$$
\D{\la}{y_1^i}G_1{}^i{}_k
=
gQ_1
\Biggl[
\fr{\fr12g\fr{v_{1k}}{q_1}A_2+h^2v_{2k}} {\sqrt{B_1}\,\sqrt{B_2}}
-\fr{\la}{B_1}
\lf(1+\fr12g\fr{b_1}{q_1}\rg)
v_{1k}
\Biggr]
$$
\ses
$$
+
gP_{1k}
\Biggl[
\fr{(A_1-b_1)A_2+h^2r_{ij}y_1^iy_2^j  -h^2b_1b_2} {\sqrt{B_1}\,\sqrt{B_2}}
-\fr{\la}{B_1}
(q_1^2+\fr12gb_1q_1)
\Biggr]
$$
\ses\ses
$$
-g\fr{v_{1k}}{q_1}
\fr{h^2} {\sqrt{B_1}\,\sqrt{B_2}}
f_1^iu_{2i}
-g\fr{v_{1k}}{q_1}
\Biggl[
\fr{\lf(1-\fr12g\fr{b_1}{q_1}\rg)A_2-h^2b_2} {\sqrt{B_1}\,\sqrt{B_2}}
-\fr{g\la}{2B_1}
\Bigl(q_1-\fr{b_1}{q_1}b_1\Bigr)
\biggr]
\Biggr]
f_1^ib_i
$$
\ses\ses
$$
+g
\Biggl[
\fr{\fr12g\lf(b_2+\fr12gq_2\rg)}    {\sqrt{B_1}\,\sqrt{B_2}}
-\fr{\la }{B_1}
(q_1+\fr12gb_1)
\Biggr]
f_{1k}
+gq_1
\fr{1-\fr14g^2} {\sqrt{B_1}\,\sqrt{B_2}}
f_{2k}
$$
\ses\ses
\be
-g
\Biggl[
\fr{(q_1-\fr12gb_{1})A_2-h^2b_2q_1} {\sqrt{B_1}\,\sqrt{B_2}}
-\fr{g\la}{2B_1}
(q_1^2+b_1^2)
\Biggr]
f^i{}_kb_i
+\De.
\ee
\ses
{
\nin
Finally,
\ses\\
$$
\D{\la}{y_1^i}G_1{}^i{}_k
=
gQ_1
\Biggl[
\fr{\fr12g\fr{v_{1k}}{q_1}A_2+h^2v_{2k}} {\sqrt{B_1}\,\sqrt{B_2}}
-\fr{\la}{B_1}
\lf(1+\fr12g\fr{b_1}{q_1}\rg)
v_{1k}
\Biggr]
$$
\ses
$$
+
gb_1P_{1k}
\Biggl[
-\fr{A_2 +h^2b_2} {\sqrt{B_1}\,\sqrt{B_2}}
+\fr{\la}{B_1}
A_1
\Biggr]
-g\fr{v_{1k}}{q_1}
\fr{h^2} {\sqrt{B_1}\,\sqrt{B_2}}
f_1^iu_{2i}
$$
\ses\ses
$$
-g\fr{v_{1k}}{q_1}
\Biggl[
\fr{\lf(1-\fr12g\fr{b_1}{q_1}\rg)A_2-h^2b_2} {\sqrt{B_1}\,\sqrt{B_2}}
-\fr{\la}{B_1}
\fr12g\Bigl(q_1-\fr{b_1}{q_1}b_1\Bigr)
\Biggr]
f_1^ib_i
$$
\ses\ses
\be
+g
\Biggl[
\fr{q_2+\fr12gb_2}    {\sqrt{B_1}\,\sqrt{B_2}}
-\fr{\la }{B_1}
(q_1+\fr12gb_1)
\Biggr]
f_{1k}
-g
\Biggl[
\fr{(q_1-\fr12gb_{1})A_2-h^2b_2q_1} {\sqrt{B_1}\,\sqrt{B_2}}
-\fr{\la}{B_1}
\fr12g(q_1^2+b_1^2)
\Biggr]
f^i{}_kb_i
+\De.
\ee
\ses

{

\ses\ses

\def\bibit[#1]#2\par{\rm\noindent\parskip1pt
                     \parbox[t]{.05\textwidth}{\mbox{}\hfill[#1]}\hfill
                     \parbox[t]{.925\textwidth}{\baselineskip11pt#2}\par}

\nin {\bf  REFERENCES}

\ses\ses



\bibit[1] H. Rund: \it The Differential Geometry of Finsler  Spaces, \rm Springer, Berlin 1959.

\bibit[2] G.S. Asanov: \it Finsler Geometry, Relativity and Gauge
 Theories, \rm D.~Reidel Publ. Comp., Dordrecht 1985.

\bibit[3] D.~Bao, S.S. Chern, and Z. Shen: {\it  An
Introduction to Riemann-Finsler Geometry,}  Springer, N.Y., Berlin
2000.

\bibit[4] G.S. Asanov:   Finsleroid space with angle and scalar product,
\it Publ.  Math.  Debrecen \bf 67 \rm(2005), 209-252.

\bibit[5] G.S. Asanov:  Finsleroid--Finsler  space with Berwald and  Landsberg conditions,
 {\it  arXiv:math.DG}/0603472 (2006).

\bibit[6] G.S. Asanov:  Finsleroid--Finsler  space and spray   coefficients, {\it  arXiv:math.DG}/0604526 (2006).

\bibit[7] G.S. Asanov:  Finsleroid--Finsler  spaces of positive--definite and  relativistic types,
\it Rep. Math. Phys. \bf 58 \rm(2006), 275--300.

\bibit[8]   G.S. Asanov:  Finsleroid--Finsler space and geodesic spray  coefficients,
\it Publ. Math.  Debrecen \bf 70 \rm(2006) (to appear).

{

\end {document}